\batchmode
\documentclass[11pt]{article}

\usepackage{epsfig}
\usepackage{graphicx}
\usepackage{color}

\newtheorem{theorem}{Theorem}
\newtheorem{lemma}{Lemma}
\newtheorem{corollary}{Corollary}

\newcommand{\be}{\begin{equation}}
\newcommand{\ee}{\end{equation}}
\newcommand{\bea}{\begin{eqnarray}}
\newcommand{\eea}{\end{eqnarray}}
\newcommand{\beas}{\begin{eqnarray*}}
\newcommand{\eeas}{\end{eqnarray*}}
\newcommand{\ba}{\begin{array}}
\newcommand{\ea}{\end{array}}

\definecolor{armygreen}{rgb}{0.29, 0.33, 0.13}

\newcommand{\real}{\mbox{$\mathrm{I\!R}$}}

\newcommand{\bfI}{\ensuremath{\mathbf{I}}}

\def\XXint#1#2#3{{\setbox0=\hbox{$#1{#2#3}{\int}$}
     \vcenter{\hbox{$#2#3$}}\kern-.5\wd0}}

\newcommand{\mcD}{\ensuremath{\mathcal{D}}}

\newcommand{\mcF}{\ensuremath{\mathcal{F}}}

\newcommand{\mcL}{\ensuremath{\mathcal{L}}}

\newcommand{\mrI}{\ensuremath{\mathrm{I}}}

\def\qed{\hbox{\vrule width 6pt height 6pt depth 0pt}}

\usepackage{amsmath}
\usepackage{amssymb}

\title{Existence and Regularity of solutions to 1-D Fractional Order Diffusion Equations} 
\author{Lueling~Jia\thanks{Applied and Computational Mathematics Division,
	  Beijing Computational Science Research Center, Beijing, China.
	  email: {\tt lljia@csrc.ac.cn}.}
       \and
         Huanzhen~Chen\thanks{School of Mathematics and Statistics,
	  Shandong Normal University, Jinan, China.
	  email: {\tt chhzh@sdnu.edu.cn}.}
	\and 
	V.J.~Ervin\thanks{Department of Mathematical Sciences,
	  Clemson University, Clemson, South Carolina 29634-0975, USA.
	  email: {\tt vjervin@clemson.edu}. Part of this work was undertaken while V.J.~Ervin was a visitor at the 
           School of Mathematics and Statistics,
	  Shandong Normal University, Jinan, China.}
	}

\date{\today}

\begin{document}
\maketitle

\begin{abstract}
In this article we investigate the existence and regularity of 1-d 
steady state fractional order diffusion equations. Two models are investigated:
the Riemann-Liouville fractional diffusion equation, and the Riemann-Liouville-Caputo
fractional diffusion equation. For these models we explicitly show how the regularity of 
the solution depends upon the right hand side (rhs) function. 
We also establish
for which Dirichlet and Neumann boundary conditions the models are well posed.
\end{abstract}

\textbf{Key words}.  Fractional diffusion equation, existence, regularity, spectral method

\textbf{AMS Mathematics subject classifications}. 35R11, 35R25, 65N35 

\setcounter{equation}{0}
\setcounter{figure}{0}
\setcounter{table}{0}
\setcounter{theorem}{0}
\setcounter{lemma}{0}
\setcounter{corollary}{0}
\section{Introduction}
 \label{sec_intro}
In recent years nonlocal models have been proposed to model a number of phenomena whose
behavior differ significantly from that predicted by usual local models, i.e., integer order
differential equations. Several areas where nonlocal models have been used include
contaminant transport in ground water flow \cite{ben001},
viscoelasticity \cite{mai971}, image processing \cite{bua101, gat151},
turbulent flow \cite{mai971, shl871}, and chaotic dynamics \cite{zas931}.

Two nonlocal approaches that are currently being investigated as models for 
anomalous diffusion are fractional differential equations \cite{pod991, kil061} and
equations involving the fractional Laplacian \cite{poz161}. (For recent results
on the regularity of the solution to equations involving the fractional Laplacian
see \cite{aco171, aco172}.) The focus of this article is on the regularity of the solution
to fractional diffusion equations. Two such models that have appeared in the literature,
which we denote by $\mbox{}_{RLC}\mcD_{r}^{\alpha} \cdot$ \cite{erv162, erv061}, and
$\mbox{}_{RL}\mcD_{r}^{\alpha} \cdot$~\cite{sch011}, are defined by
\begin{align}
 \mbox{}_{RLC}\mcD_{r}^{\alpha} u(x)  &:= \,  
   - D \left( r D^{-(2 - \alpha)} \, + \, (1 - r) D^{-(2 - \alpha)*}  \right) D u(x) \ = \ f(x) \, , \ \ 0 < x < 1 \, , \
   \label{DefDRLC}  \\
   \mbox{ and } \quad \mbox{}_{RL}\mcD_{r}^{\alpha} u(x) &:= \,  
   - D^{2} \left( r D^{-(2 - \alpha)} \, + \, (1 - r)  D^{-(2 - \alpha)*}  \right)  u(x) \ = \ f(x) \, , \ \ 0 < x < 1 \, ,   
\label{DefDRL}
\end{align}   
where $D$ denotes the usual differential operator, and $D^{- \beta}$ and $D^{- \beta *}$ denote
the left and right fractional integral operators, respectively (defined in Section \ref{sec_not}).
We refer to $\mbox{}_{RLC}\mcD_{r}^{\alpha}$ as the Riemann-Liouville-Caputo fractional
differential operator, and $ \mbox{}_{RL}\mcD_{r}^{\alpha}$ as the Riemann-Liouville fractional
differential operator

In \cite{erv061} a variational solution, $u \in H^{\alpha/2}_{0}(0 , 1)$, to \eqref{DefDRLC} and
\eqref{DefDRL}, subject to $u(0) = u(1) = 0$, and $f \in H^{-\alpha/2}(0 , 1)$ was established.
A detailed analysis of the existence and regularity of solutions to  \eqref{DefDRLC} and
\eqref{DefDRL} for $ r = 1$ (i.e., one sided fractional diffusion equations) was given by
Jin, Lazarov, et al. in \cite{jin151}. Recently in \cite{wan171} Wang and Yang investigated the
well posedness of solutions to  \eqref{DefDRLC} and
\eqref{DefDRL} for $ r = 1$ subject to three different Neumann boundary conditions. They 
showed that for at least one of the boundary conditions that the modeling equations were
ill posed. A physical interpretation of absorbing and reflecting boundary conditions for
\eqref{DefDRL} for $ r = 1$ was recently presented by Baeumer, Kov\'{a}cs, et al. in \cite{bae181}.
Also, for $r = 1$ the existence and regularity 
of solutions to \eqref{DefDRLC} having a variable diffusion coefficient
was analyzed by Yang, Chen and Wang in \cite{yan181}.

In this paper we investigate the existence and regularity of the solutions to 
\eqref{DefDRLC} and \eqref{DefDRL}, subject to various boundary conditions. A goal of this
investigations is to provide engineers and scientists insight into determining which equation
may more appropriately model their problem of interest. 

For clarity in our discussion we say that $g(x), \, x \in (0 , 1)$ is \textit{algebraically regular} if
$g(x) \sim C \, x^{a}$, as $x \rightarrow 0$ for $0 < a < 1$ or 
$g(x) \sim C \, (1 - x)^{b}$, as $x \rightarrow 1$ for $0 < b < 1$, and
\textit{algebraically singular} if
$g(x) \sim C \, x^{- a}$, as $x \rightarrow 0$ for $0 < a < 1$ or 
$g(x) \sim C \, (1 - x)^{- b}$, as $x \rightarrow 1$ for $0 < b < 1$.

From an elementary course in differential equations we have that the general solution to the linear
differential equation $\mcL u \, = \, f$ can be expressed as $u \, = \, u_{homog} \, + \, u_{ps}$, 
where $u_{homog}$ satisfies the associated homogeneous differential equation (i.e., $u_{homog} \in
ker( \mcL )$), and $u_{ps}$ is a particular solution. The regularity of the solution $u$ depends on 
two factors: (i) the operator $\mcL$, and (ii) the RHS function $f$. The regularity of $u_{homog}$
is solely determined by the operator $\mcL$. The regularity of $u_{ps}$ depends on $f$, and
also on the operator $\mcL$. For \eqref{DefDRLC} we have that 
$ker( \mbox{}_{RLC}\mcD_{r}^{\alpha} )$ is algebraically regular (Section \ref{sec_exRLC}), whereas for
\eqref{DefDRL} we have that $ker( \mbox{}_{RL}\mcD_{r}^{\alpha} )$ is algebraically 
singular (Section \ref{sec_exRL}).

Equation \eqref{DefDRL} represents the steady-state fractional diffusion equation derived in \cite{sch011},
assuming a heavy tail random walk process. Equation \eqref{DefDRLC} and \eqref{DefDRL} only differ
in the location of one of the derivative operators: either before the fractional integral terms or
after them. A physical interpretation of the difference between the two equations can be
obtained by considering the 1-D heat equation, modeling the cross sectional temperature along a
bar that is insulated along its lateral surface \cite{erv160}.
\be
\frac{\partial}{\partial t}u(x, t) \ - \ \frac{\partial}{\partial x}q(x , t) \ = \ f(x , t) \, , \quad 0 < x < 1, \ \ t > 0 \, .
\label{deheat}
\ee
In \eqref{deheat} $u(x , t)$, $q(x , t)$, and $f(x , t)$ represent the temperature (synonymous with energy),
energy flux, and an energy source density, respectively,  at cross section $x$ at time $t$. Corresponding to
\eqref{DefDRLC},
\[
q(x , t) \ = \ \frac{1}{\Gamma(2 - \alpha)} \left( r \, \int_{0}^{x} \frac{1}{(x - s)^{\alpha - 1}} 
(- \frac{\partial u(s, t)}{\partial x} ) \,  ds \ - \ 
(1 - r)\, \int_{x}^{1} \frac{1}{(s - x)^{\alpha - 1}} 
 \, \frac{\partial u(s, t)}{\partial x}  \,  ds \right) \, .
\]
The first term in the parenthesis on the right hand side implies that the if the local temperature around $s$
is not constant (i.e., $ \frac{\partial u(s, t)}{\partial x} \ne 0 $) then energy flows from this point. The contribution
of this flow of energy to a point a distant $(x - s)$ units away is given by $\frac{1}{\Gamma(2 - \alpha)} \, 
\frac{1}{(x - s)^{\alpha - 1}} 
(- \frac{\partial u(s, t)}{\partial x} ) \,  \Delta s$. A similar interpretation applies to the second term in the
parenthesis. Hence, in \eqref{DefDRLC} the flux at a point is the weighted sum of local energy variations
along the bar, which may be interpreted as a nonlocal version of Fick's Law.

With equation \eqref{DefDRL}, let
\[
   E(x , t) \ = \ \frac{1}{\Gamma(2 - \alpha)} \left( r \, \int_{0}^{x} \frac{1}{(x - s)^{\alpha - 1}} 
\, u(s, t)  \,  ds \ + \ 
(1 - r)\, \int_{x}^{1} \frac{1}{(s - x)^{\alpha - 1}} 
 \,   u(s, t)  \,  ds \right) \, .
\]
$E(x , t)$ may be interpreted as the weighted sum of local energy distributed throughout the bar, with
the energy (temperature) at $s$ contributing an amount 
$\frac{1}{\Gamma(2 - \alpha)} \, \frac{1}{(x - s)^{\alpha - 1}}  u(s, t)  \,  \Delta s$. Then, as
$q(x , t) \ = \ - \, \frac{\partial}{\partial x} E(x , t)$, the flux at $x$ is due to the variation in the 
weighted energy at $x$. Note that each point $s$ contributes to the weighted energy at $x$, $E(x , t)$,
corresponding to a random walk process as derived in \cite{sch011}. So, in \eqref{DefDRL} there
is an underlying energy flow occurring throughout the bar, however there is only a resulting flux at
$x$ if there is a local imbalance in this weighted energy at $x$.

Following in the next section we introduce notation and several key lemmas we use in the analysis of 
the solutions to \eqref{DefDRLC} and \eqref{DefDRL}. In Section \ref{sec_exRLC} we present the existence
and regularity results for the solution of \eqref{DefDRLC}, subject to various boundary conditions. 
A shift theorem for \eqref{DefDRLC} is investigated in Section \ref{ssecRegu}. 
The
analysis of the solution to \eqref{DefDRL}, subject to various boundary conditions, is presented in 
Section \ref{sec_exRL}. A summary of the difference in the solutions of \eqref{DefDRLC} and \eqref{DefDRL}
is given in the Conclusions.
Proofs of a number of the results used in Sections \ref{sec_exRLC} and
\ref{sec_exRL} are given in the appendix.

\setcounter{equation}{0}
\setcounter{figure}{0}
\setcounter{table}{0}
\setcounter{theorem}{0}
\setcounter{lemma}{0}
\setcounter{corollary}{0}
\section{Notation and Properties}
\label{sec_not}
For $u$ a function defined on $(a , b)$, and $\sigma > 0$,
we have that the left and right fractional integral operators are defined as: \\
\underline{Left Fractional Integral Operator}:
 $ \mbox{}_{a}D_{x}^{-\sigma}u(x) \, := \, \frac{1}{\Gamma(\sigma)} \int_{a}^{x} (x - s)^{\sigma - 1} \, u(s) \, ds \, . $ 
 \vspace{0.5em} \\
\underline{Right Fractional Integral Operator}:
 $ \mbox{}_{x}D_{b}^{-\sigma}u(x)  \, := \, \frac{1}{\Gamma(\sigma)} \int_{x}^{b} (s - x)^{\sigma - 1} \, u(s) \, ds \, . $ 
 \vspace{0.5em} \\
Then, for $\mu > 0$, $n$ the smallest integer greater than $\mu$ $(n-1 \le \mu < n)$,  
$\sigma = n - \mu$, and $D$ the derivative operator,
the left and right  Riemann-Liouville fractional differential operators are defined as:  \\
\underline{Left Riemann-Liouville Fractional Differential Operator of order $\mu$}: \\ 
 \mbox{} \hfill $ \mbox{}_{a}^{RL}D_{x}^{\mu} u(x) \, := \, D^{n} \mbox{}_{a}D_{x}^{-\sigma}u(x) \ = \ 
 \frac{1}{\Gamma(\sigma)} \frac{d^{n}}{dx^{n}} \int_{a}^{x} (x - s)^{\sigma - 1} \, u(s) \, ds \, . $ 
 \vspace{0.5em} \\
\underline{Right Riemann-Liouville Fractional Differential Operator of order $\mu$}: \\
 \mbox{} \hfill  $ \mbox{}_{x}^{RL}D_{b}^{\mu} \, u(x) := \, (- D)^{n} \mbox{}_{x}D_{b}^{-\sigma}u(x) \ = \ 
 \frac{(-1)^{n}}{\Gamma(\sigma)} \frac{d^{n}}{dx^{n}} \int_{x}^{b} (s - x)^{\sigma - 1} \, u(s) \, ds \, . $ 
 
The Riemann-Liouville and Caputo fractional differential operators differ in the location of the derivative operator. \\
\underline{Left Caputo Fractional Differential Operator of order $\mu$}: \\
 \mbox{} \hfill $ \mbox{}_{a}^{C}D_{x}^{\mu} u(x) \, := \,  \mbox{}_{a}D_{x}^{-\sigma}\, D^{n} u(x) \ = \ 
 \frac{1}{\Gamma(\sigma)} \int_{a}^{x} (x - s)^{\sigma - 1} \, \frac{d^{n}}{ds^{n}} u(s) \, ds \, . $ 
 \vspace{0.5em} \\
\underline{Right Caputo Fractional Differential Operator of order $\mu$}: \\
 \mbox{} \hfill  $ \mbox{}_{x}^{C}D_{b}^{\mu} \, u(x) := \, (-1)^{n} \mbox{}_{x}D_{b}^{-\sigma}\, D^{n} u(x) \ = \ 
 \frac{(-1)^{n}}{\Gamma(\sigma)} \int_{x}^{b} (s - x)^{\sigma - 1} \, \frac{d^{n}}{ds^{n}} u(s) \, ds \, . $ 

As our interest is in the solution of fractional diffusion equations on a bounded, connected subinterval of $\real$, 
without loss
of generality we restrict our attention to the unit interval $\mrI \, := \, (0 , 1)$.

For ease of notation, we use 
\[
   D^{-\sigma} \, := \,  \mbox{}_{0}D_{x}^{-\sigma} \, ,  \ \mbox{ and } \ 
    D^{-\sigma *} \, := \,  \mbox{}_{x}D_{1}^{-\sigma} \, .
\] 

Let,
\be
 \bfI_{r}^{s}u(x)  \ := \ r  D^{-s} u(x) \, + \, (1 - r)  D^{-s*} u(x) \, .
\label{defbfI}
\ee

Then,
\[
    \mbox{}_{RLC}\mcD_{r}^{\alpha} u(x) \ = \ - D \bfI_{r}^{2 - \alpha} D u(x) \, , \quad \mbox{ and }
      \mbox{}_{RL}\mcD_{r}^{\alpha} u(x) \ = \ - D^{2} \bfI_{r}^{2 - \alpha} u(x) \, .
\]

For the $RLC$ fractional diffusion equation, the \textit{flux} is given by 
$\mbox{}_{RLC}\mcF u(x)  \ = \ -   \bfI_{r}^{2 - \alpha} D u(x)$ and for 
the $RL$ fractional diffusion equation, the \textit{flux} is given by 
$\mbox{}_{RL}\mcF u(x)  \ = \ -   D \bfI_{r}^{2 - \alpha} u(x)$.

Jacobi polynomial play a key role in the approximation schemes. We briefly review their definition and
properties central to the method \cite{abr641, sze751}. 

\textbf{Usual Jacobi Polynomials, $P_{n}^{(\alpha , \beta)}(x)$, on $(-1 \, , \, 1)$}.   \\    
\underline{Definition}: $ P_{n}^{(\alpha , \beta)}(x) \ := \ 
\sum_{m = 0}^{n} \, p_{n , m} \, (x - 1)^{(n - m)} (x + 1)^{m}$, where
\begin{equation}
       p_{n , m} \ := \ \frac{1}{2^{n}} \, \left( \begin{array}{c}
                                                              n + \alpha \\
                                                              m  \end{array} \right) \,
                                                    \left( \begin{array}{c}
                                                              n + \beta \\
                                                              n - m  \end{array} \right) \, .
\label{spm21}
\end{equation}
\underline{Orthogonality}:    
\begin{align}
 & \int_{-1}^{1} (1 - x)^{\alpha} (1 + x)^{\beta} \, P_{j}^{(\alpha , \beta)}(x) \, P_{k}^{(\alpha , \beta)}(x)  \, dx 
 \ = \
   \left\{ \begin{array}{ll} 
   0 , & k \ne j  \\
   |\| P_{j}^{(\alpha , \beta)} |\|^{2}
   \, , & k = j  
    \end{array} \right.  \, ,  \nonumber \\
& \quad \quad \mbox{where } \  \ |\| P_{j}^{(\alpha , \beta)} |\| \ = \
 \left( \frac{2^{(\alpha + \beta + 1)}}{(2j \, + \, \alpha \, + \, \beta \, + 1)} 
   \frac{\Gamma(j + \alpha + 1) \, \Gamma(j + \beta + 1)}{\Gamma(j + 1) \, \Gamma(j + \alpha + \beta + 1)}
   \right)^{1/2} \, .
  \label{spm22}
\end{align}                                                    

In order to transform the domain of the family of Jacobi polynomials to $[0 , 1]$, let $x \rightarrow 2t - 1$ and 
introduce $G_{n}^{(\alpha , \beta)}(t) \, = \, P_{n}^{(\alpha , \beta)}( x(t) )$. From \eqref{spm22},
\begin{align}
 \int_{-1}^{1} (1 - x)^{\alpha} (1 + x)^{\beta} \, P_{j}^{(\alpha , \beta)}(x) \, P_{k}^{(\alpha , \beta)}(x)  \, dx 
 &= \
 \int_{t = 0}^{1} 2^{\alpha} \, (1 - t)^{\alpha} \, 2^{\beta} \, t^{\beta} \, P_{j}^{(\alpha , \beta)}(2t - 1) \, P_{k}^{(\alpha , \beta)}(2t - 1)  \, 2 \,  dt
 \nonumber \\
  &= \
2^{\alpha + \beta + 1} \int_{t = 0}^{1}   (1 - t)^{\alpha}  \, t^{\beta} \, G_{j}^{(\alpha , \beta)}(t) \, G_{k}^{(\alpha , \beta)}(t)  \,  dt
\nonumber \\
&= \
   \left\{ \begin{array}{ll} 
   0 , & k \ne j \, , \\
  2^{\alpha + \beta + 1} \, |\| G_{j}^{(\alpha , \beta)} |\|^{2}
   \, , & k = j  
   \, . \end{array} \right.    \nonumber \\
 \quad \quad \mbox{where } \  \ |\| G_{j}^{(\alpha , \beta)} |\| &= \
 \left( \frac{1}{(2j \, + \, \alpha \, + \, \beta \, + 1)} 
   \frac{\Gamma(j + \alpha + 1) \, \Gamma(j + \beta + 1)}{\Gamma(j + 1) \, \Gamma(j + \alpha + \beta + 1)}
   \right)^{1/2} \, .  \label{spm22g} 
\end{align}                                                    

\begin{equation}
\mbox{Note that } \quad  |\| G_{j}^{(\alpha , \beta)} |\| \ = \ |\| G_{j}^{(\beta , \alpha)} |\| \, ,
\label{nmeqG}
\end{equation}
and from \cite{abr641, sze751}
\begin{equation}
G_{j}^{(\alpha , \beta)}(0) \ = \ (-1)^{j} \, \frac{\Gamma(j + \beta + 1)}{\Gamma(j + 1) \, \Gamma(\beta + 1)} \, .
\label{eqG0}
\end{equation}

 From \cite[equation (2.19)]{mao161} we have that
\begin{equation}
   \frac{d^{k}}{dx^{k}} P_{n}^{(\alpha , \beta)}(x) \ = \ 
   \frac{\Gamma(n + k + \alpha + \beta + 1)}{2^{k} \, \Gamma(n + \alpha + \beta + 1)} P_{n - k}^{(\alpha + k \, , \, \beta + k)}(x) \, .
   \label{derP}
\end{equation}   
Hence,
\begin{align}
\frac{d^{k}}{dt^{k}} G_{n}^{(\alpha , \beta)}(t) 
  &= \ \frac{\Gamma(n + k + \alpha + \beta + 1)}{  \Gamma(n + \alpha + \beta + 1)} 
  G_{n - k}^{(\alpha + k \, , \, \beta + k)}(t)  \, .  \label{eqC4}
\end{align}   

Also, from \cite[equation (2.15)]{mao161}, 
\begin{equation}
\frac{d^{k}}{dx^{k}} \left\{ (1 - x)^{\alpha + k} \, (1 + x)^{\beta + k} \, P_{n - k}^{(\alpha + k \, , \, \beta + k)}(x) \right\}
 \ = \ 
 \frac{(-1)^{k} \, 2^{k} \, n!}{(n - k)!} \, (1 - x)^{\alpha} \, (1 + x)^{\beta} \, P_{n}^{(\alpha \, , \, \beta)}(x) \, , \
 n \ge k \ge 0 \, ,
 \label{eqB0}
\end{equation}
from which it follows that
\begin{equation}
 \frac{d^{k}}{dt^{k}} \left\{  \ (1 \, - \, t)^{\alpha + k} \,  t^{\beta + k} \,
 G_{n - k}^{(\alpha + k \, , \, \beta + k)}(t) \right\} 
 \ = \ 
 \frac{(-1)^{k} \,  n!}{(n - k)!} \,  (1 \, - \, t)^{\alpha} \,  t^{\beta} \,
 G_{n}^{(\alpha \, , \, \beta)}(t) \, . 
  \label{eqC2}
 \end{equation}
 
For compactness of notation we introduce
\begin{equation}
 \rho^{(\alpha , \beta)} \, = \, \rho^{(\alpha , \beta)}(x) \, := \, (1 - x)^{\alpha} \, x^{\beta} \, .
 \label{defrho}
\end{equation} 

 We use $y_{n} \sim n^{p}$ to denote that there exists constants $c$ and $C > 0$ such that, as 
 $n \rightarrow \infty$,  \linebreak[4]
 $c \, n^{p} \le | y_{n} | \le C \, n^{p}$. Also, 
from Stirling's formula we have that
\begin{equation}
\lim_{n \rightarrow \infty} \, \frac{\Gamma(n + \sigma)}{\Gamma(n) \, n^{\sigma}}
\ = \ 1 \, , \mbox{ for } \sigma \in \mathbb{R}.  
 \label{eqStrf}
\end{equation}

\textbf{Function spaces, $L_{\omega}^{2}(\mrI)$ and $H^{l}_{\rho^{(a \, , \, b)} \,  , \, A}(\mrI)$}. \\
The weighted $L^{2}((\mrI)$ spaces are appropriate for analyzing the convergence of the spectral type methods 
presented below. For $\omega(x) > 0, \ x \in (0 , 1)$, let 
\[
L_{\omega}^{2}(0 , 1) \, := \, \{ f(x) \, : \, \int_{0}^{1} \omega(x) \, f(x)^{2} \, dx \ < \ \infty \} \, .
\]
Associated with $L_{\omega}^{2}(0 , 1)$ is the inner product, $\langle \cdot , \cdot \rangle_{\omega}$, and
norm, $\| \cdot \|_{\omega}$, defined by
\begin{align*}
\langle f \,  , \, g \rangle_{\omega} &:= \ \int_{0}^{1} \omega(x) \, f(x) \, g(x) \, dx \, , \quad \mbox{and} \\
 \| f \|_{\omega} &:= \ \left( \langle f \,  , \, f \rangle_{\omega} \right)^{1/2} \, .
\end{align*}

Following \cite{guo041}, we introduce weighted Sobolev spaces $H^{l}_{\rho^{(a \, , \, b)} \,  , \, A}(\mrI)$ defined by
\be
  H^{l}_{\rho^{(a \, , \, b)} \,  , \, A}(\mrI) \ := \
  \left\{ v \, | \, v \mbox{ is measurable and } \| v \|_{l , \, \rho^{(a \, , \, b)} , \, A} < \infty \right\} \, , \ \ l \in \mathbb{N} \, ,
\label{defwHr}
\ee
with associated norm and semi-norm
\[
   \| v \|_{l , \, \rho^{(a \, , \, b)} , \, A}  \ := \ \left( \sum_{j = 0}^{l} 
   \| D^{j} v \|_{\rho^{(a + j \, , b + j)}}^{2} \right)^{1/2} \, , \ \ \
    | v |_{l , \, \rho^{(a \, , \, b)} , \, A}  \ := \  \| D^{l} v \|_{\rho^{(a + l \, , b + l)}} \, .
\]

Throughout the manuscript we assume that $\alpha$, $\beta$, and $r$ satisfy a fixed relationship. Additionally
a constant defined by $\alpha$ and $\beta$ occurs sufficiently often that we denote it by $c_{*}^{*}$.  We 
refer to these relationships as \textbf{Condition A}.

\underline{Definition}: \textbf{Condition A} \\
The parameters $\alpha$, $\beta$, and $r$ and constant $c_{*}^{*}$ satisfy:
$1 < \alpha < 2$, $\alpha - 1 \, \le  \, \beta \, , \, \alpha - \beta \, \le  \, 1$, $0 \le r \le 1$ 
\begin{equation}
  c_{*}^{*} \ = \ \frac{\sin(\pi \alpha)}{\sin(\pi (\alpha - \beta)) \, + \, \sin(\pi \beta)} \, ,  \label{defcss}
\end{equation}  
where $\beta$ is determined by
\begin{equation}
  r \ = \ \frac{\sin( \pi \, \beta)}{\sin( \pi ( \alpha - \beta)) \, + \,  \sin( \pi \, \beta)} \, . \label{propker0} 
\end{equation}

The following three lemmas are useful in determining the solutions of \eqref{DefDRLC} and \eqref{DefDRL}.

\begin{lemma}   \label{ker1r}
Under \textbf{Condition A}, 
\begin{equation}
  ker(D \bfI_{r}^{2 - \alpha}) \ = \ span \{ \rho^{(\alpha - \beta - 1 \, , \,   \beta - 1)}(x) \}  \, . 
  \label{ytre1}
  \end{equation}
Additionally, 
\begin{align}
 D \bfI_{r}^{2 - \alpha} (x \, \rho^{(\alpha - \beta - 1 \, , \,   \beta - 1)}(x) ) &= \  - \, c_{*}^{*} \, \Gamma(\alpha)   
 \ = \ \mu_{-1} \, G_{0}^{(\delta , \gamma)}(x) \, ,   \label{chr11} \\
 \mbox{and } \ \   D \bfI_{r}^{2 - \alpha} ((1 - x) \, \rho^{(\alpha - \beta - 1 \, , \,   \beta - 1)}(x)) 
 &= \, - \mu_{-1} \, G_{0}^{(\delta , \gamma)}(x) \, ,   
 \quad \mbox{where } 
   \mu_{-1} \,  := \,  - \, c_{*}^{*} \, \Gamma(\alpha)  \, .
  \label{chr12} 
\end{align}  
\end{lemma}
\textbf{Proof}: The proof of \eqref{ytre1} is given in  \cite{erv162}. Properties \eqref{chr11} and \eqref{chr12} follow from
Lemma \ref{spec_poly} (in the appendix), and that  $G_{0}^{(\delta , \gamma)}(x) \, = \, 1$. \\
\mbox{ } \hfill \qed

\begin{lemma} \label{lmacherry2}
Under \textbf{Condition A}, 
for $n \, = \, 0, 1, 2, \ldots$
\begin{align}
 \bfI_{r}^{2 - \alpha} (1 - x)^{\alpha - \beta - 1} \, x^{\beta - 1} \,  G_{n}^{(\alpha - \beta - 1 \, , \, \beta - 1)}(x) \, 
 &= \  \sigma_{n} \, G_{n}^{(\beta - 1 \, , \, \alpha - \beta - 1)}(x) \, ,    \label{chr111} \\
 \quad \mbox{where } 
   \sigma_{n} \, &:= \,  - \, c_{*}^{*} \, \frac{\Gamma(n + \alpha - 1)}{\Gamma(n + 1)} \, .
  \label{chr113} 
\end{align}   
\end{lemma}
The proof of this lemma is given in the appendix.

\begin{lemma} \cite{mao181} \label{lmaegp}
Under \textbf{Condition A}, 
for $n \, = \, 0, 1, 2, \ldots$
\begin{align}
 D \bfI_{r}^{2 - \alpha} (1 - x)^{\alpha - \beta} x^{\beta} \, G_{n}^{(\alpha - \beta \, , \, \beta)}(x) 
 &= \ \mu_{n}  \, G_{n + 1}^{(\beta - 1 \, , \, \alpha - \beta - 1)}(x) \, ,  \label{NactG}  \\
 \mbox{where }  \quad
   \mu_{n} &= \ \, c_{*}^{*} \, \, \frac{\Gamma(n + \alpha)}{\Gamma(n + 1)} \, .
 \label{defln}
\end{align}
\end{lemma}
\mbox{ } \hfill \qed 
%

\setcounter{equation}{0}
\setcounter{figure}{0}
\setcounter{table}{0}
\setcounter{theorem}{0}
\setcounter{lemma}{0}
\setcounter{corollary}{0}
\section{Existence and Regularity of the RLC Fractional Diffusion Model}
\label{sec_exRLC}
In this section we investigate the existence and regularity of the steady state RLC fractional diffusion model subject to
various boundary conditions.

\subsection{Dirichlet Boundary Conditions}
\label{ssec_RLCdbc}
From \cite{erv162} we have that
\begin{align}
 ker( \mbox{}_{RLC}\mcD_{r}^{\alpha} ) &= \ span \left\{ 1 \ ,  \ \int_{0}^{x} (1 - s)^{\alpha - \beta - 1} \, s^{\beta - 1} \, ds 
 \right\}  \nonumber \\
\ = \ span \left\{ k_{0}(x) \,  \right. &:= \,  \left. \int_{x}^{1} (1 - s)^{\alpha - \beta - 1} \, s^{\beta - 1} \, ds  \ , 
\ k_{1}(x) \, := \, \int_{0}^{x} (1 - s)^{\alpha - \beta - 1} \, s^{\beta - 1} \, ds   \right\} \, .
     \label{kRLC1}
\end{align}

The singular endpoint behavior of the kernel at both endpoints, i.e., $(1 - x)^{\alpha - \beta}$ and $x^{\beta}$,
is more apparent using the basis $k_{0}(x)$ and $k_{1}(x)$.
   
With $C_{1}$ and $C_{2}$ appropriately chosen, the change of variable $\tilde{u}(x) \ = \ u(x) \, + \, C_{1} k_{0}(x) \,
+ \, C_{2} k_{1}(x)$ transform the problem
\be
    \mbox{}_{RLC}\mcD_{r}^{\alpha} \tilde{u}(x) \ = \ f(x) \, , \ \ 0 < x < 1 \, , \ \ \ 
    \mbox{subject to } \tilde{u}(0) = A \, , \ \ \tilde{u}(1) = B \, ,
   \label{rlcde1}
\ee
to the problem
\be
    \mbox{}_{RLC}\mcD_{r}^{\alpha} u(x) \ = \ f(x) \, , \ \ 0 < x < 1 \, , \ \ \ 
    \mbox{subject to } u(0) = 0 \, , \ \ u(1) = 0 \, .
\label{RLChDbc}
\ee

Note that $f(x) \in L^{2}_{\rho^{(\beta , \alpha - \beta)}}(\mrI)$ may be expressed as 
$f(x) \ = \ \sum_{i = 0}^{\infty} \frac{f_{i}}{ |\| G_{i}^{(\beta , \alpha - \beta)} |\|^{2}} \, G_{i}^{(\beta , \alpha - \beta)}(x)$, 
where $f_{i}$ is given by
\be
f_{i} \, := \, \int_{0}^{1} \,  \rho^{(\beta , \alpha - \beta)}(x) \, f(x) \, G_{i}^{(\beta , \alpha - \beta)}(x)  \, dx \, .
\label{deffis}
\ee

With $f_{i}$ defined in \eqref{deffis}, let 
\be
f_{N}(x) \ = \ \sum_{i = 0}^{N} \frac{f_{i}}{ |\| G_{i}^{(\beta , \alpha - \beta)} |\|^{2}} \, G_{i}^{(\beta , \alpha - \beta)}(x) 
 \  \mbox{ and } 
u_{N}(x) \ = \ \rho^{( \alpha - \beta , \beta )}(x) \, \sum_{i = 0}^{N} c_{i} \, G_{i}^{( \alpha - \beta , \beta )}(x) \, , 
 \label{defuNs} 
\ee
\be
\mbox{ where} \ 
\lambda_{i} \ = \ - c_{*}^{*} \,
\frac{\Gamma(i + 1 + \alpha)}{\Gamma(i + 1)}  \  \mbox{ and }
c_{i}  \ = \ \frac{1}{ \lambda_{i} \, |\| G_{i}^{( \beta , \alpha - \beta  )} |\|^{2}  }f_{i}  \, .
\label{defuNs2}
\ee

\begin{theorem}[\cite{erv162}]   \label{thmuspg}
Let $f (x) \in L^{2}_{\rho^{(\beta , \alpha - \beta)}}(\mrI)$ and $u_{N}(x)$ be as defined in \eqref{defuNs}. Then,
$u(x) \ := \ \lim_{N \rightarrow \infty} u_{N}(x) \ = \ \rho^{( \alpha - \beta , \beta )}(x) \, \sum_{j = 0}^{\infty}  c_{j} \,
G_{j}^{( \alpha - \beta , \beta )}(x) \, 
\in L^{2}_{\rho^{(-(\alpha - \beta) , -\beta)}}(\mrI)$.
In addition, $u(x)$ satisfies \eqref{RLChDbc}.
\end{theorem}
\mbox{ } \hfill \qed

The regularity of $Du$ is given by the following corollary.
\begin{corollary}  \label{corRegDu}
For $f (x) \in L^{2}_{\rho^{(\beta , \alpha - \beta)}}(\mrI)$ and $u(x)$ satisfying  \eqref{RLChDbc} we have that
$D u \in L^{2}_{\rho^{(-(\alpha - \beta) + 1 \, , \, -\beta + 1)}}(\mrI)$.
\end{corollary}
\textbf{Proof}: 
Consider 
\begin{align}
& \|  D u_{M}  \ - \ D u_{N}  \|_{\rho^{(-(\alpha - \beta) + 1 \, , \, -\beta + 1)}}^{2}  \nonumber \\
& \quad = \
 \left(\rho^{(\alpha - \beta - 1 \, , \, \beta - 1)}(x)  \, \sum_{i = N + 2}^{M + 1} c_{i - 1} \, i \, G_{i}^{(\alpha - \beta - 1 \, , \, \beta - 1)}(x) \, , \,
    \sum_{i = N + 2}^{M + 1} c_{i - 1} \, i \, G_{i}^{(\alpha - \beta - 1 \, , \, \beta - 1)}(x)  \right)   \nonumber \\
& \quad = \
    \sum_{i = N + 2}^{M + 1} c_{i - 1}^{2} \, i^{2} \, |\| G_{i}^{(\alpha - \beta - 1 \, , \, \beta - 1)}(x) |\|^{2}   \, . \nonumber
\end{align}
From \eqref{gge2}, $\frac{1}{2} \, |\| G_{i}^{(\beta - 1 \, , \, \alpha - \beta - 1)}(x) |\|^{2} \ \le \
  |\| G_{i - 1}^{(\alpha - \beta  \, , \,  \beta )}(x) |\|^{2} $, thus
\begin{align}
 \|  D u_{M}  \ - \ D u_{N}  \|_{\rho^{(-(\alpha - \beta - 1) \, , \, -(\beta - 1))}}^{2}  
& \quad \le 2 \,  
      \sum_{i = N + 2}^{M + 1}  i^{2} \, \frac{f_{i - 1}^{2}}{\lambda_{i - 1}^{2} \,  |\| G_{i - 1}^{(\beta \, , \, \alpha - \beta)}(x) |\|^{4}} \,
       |\| G_{i - 1}^{(\beta \, , \, \alpha - \beta )}(x) |\|^{2}  \nonumber \\
& \quad = \ 2 \,
      \sum_{i = N +1}^{M}  \frac{(i + 1)^{2}}{\lambda_{i}^{2}} \, \frac{f_{i}^{2}}{ |\| G_{i}^{(\beta \, , \, \alpha - \beta )}(x) |\|^{2}} \, .
      \nonumber  \\
& \quad = \ 2 \,  
 \left(\rho^{(\beta  \, , \, \alpha - \beta )}(x)  \, 
 \sum_{i = N + 1}^{M} \frac{(i + 1)^{2}}{\lambda_{i}^{2}} \, 
  \frac{f_{i}}{ |\| G_{i}^{(\beta \, , \, \alpha - \beta )}(x) |\|^{2}} \, \, G_{i}^{(\beta \, , \, \alpha - \beta )}(x) \, , \,
  \right. \nonumber \\
 & \quad \quad \quad \quad \quad \quad \quad \quad \quad \quad  
    \left. \sum_{i = N + 1}^{M}  \frac{f_{i}}{ |\| G_{i}^{(\beta \, , \, \alpha - \beta )}(x) |\|^{2}} \, \, G_{i}^{(\beta \, , \, \alpha - \beta )}(x) 
     \,  \right)   \, .  \nonumber 
\end{align}
Using Stirling's formula \eqref{eqStrf},
$  \frac{(i + 1)^{2}}{\lambda_{i}^{2}}  $ is bounded as $i \rightarrow \infty$. Hence
\begin{align}
& \|  D u_{M}  \ - \ D u_{N}  \|_{\rho^{(-(\alpha - \beta) + 1 \, , \, -\beta + 1)}}^{2}  \nonumber \\
& \quad \le C \, 
 \left(\rho^{(\beta  \, , \, \alpha - \beta )}(x)  \, 
 \sum_{i = N + 1}^{M} \, 
  \frac{f_{i}}{ |\| G_{i}^{(\beta \, , \, \alpha - \beta )}(x) |\|^{2}} \, \, G_{i}^{(\beta \, , \, \alpha - \beta )}(x)  \, , \,
   \right.  \nonumber \\
&  \quad \quad \quad \quad \quad \quad \quad \quad \quad \quad  \quad \quad \quad \quad \quad
   \left.   \sum_{i = N + 1}^{M}  \frac{f_{i}}{ |\| G_{i}^{(\beta \, , \, \alpha - \beta )}(x) |\|^{2}} \, \, G_{i}^{(\beta \, , \, \alpha - \beta )}(x) 
    \,  \right) \, \nonumber \\
& \quad = \ C    \|   f_{M}(x)  \ - \  f_{N}(x)  \|_{\rho^{(\beta  \, , \, \alpha - \beta ))}}^{2} \, .
\end{align}
 As $f \in L^{2}(\mrI)_{\rho^{(\beta \, , \, \alpha - \beta)}}$, then $\{ f_{n} \}$ is a Cauchy sequence in 
$L^{2}_{\rho^{(\beta  \, , \, \alpha - \beta )}}(\mrI)$. Thus we can conclude that 
$D u \in  L^{2}_{\rho^{( -(\alpha - \beta) + 1 \, , \, -\beta + 1)}}(\mrI)$. \\
\mbox{ } \hfill \qed

\subsection{Regularity of $u(x)$}
\label{ssecRegu}

Next we investigate if $f(x)$ is ``nicer'', i.e., a more regular function, does that increased regularity transfer over to the 
solution $u(x)$.

The following lemma is useful in helping to provide an answer to that question.

\begin{lemma} \label{thmrgqd}
For $j \in \mathbb{N}$, if $D^{j - 1}f \in L^{2}_{\rho^{(\beta + j - 1\, , \, \alpha - \beta + j - 1)}}(\mrI)$, then
$D^{j}  \frac{1}{\rho^{(\alpha - \beta  \, , \, \beta )}(x)}u(x) \,  \in L^{2}_{\rho^{(\alpha - \beta + j  \, , \, \beta + j)}}(\mrI)$.
\end{lemma}
\textbf{Proof}:
From \eqref{eqC4}, and \eqref{defuNs}
\begin{align}
D^{j} \big( \frac{1}{\rho^{(\alpha - \beta \, , \, \beta )}(x)} u_{N}(x) \big)
&= \ D^{j}  \, \sum_{i = 0}^{N} c_{i} \, G_{i}^{(\alpha - \beta \, , \, \beta)}(x) 
 \nonumber  \\
&=  \ \sum_{i = 0}^{N - j} c_{i + j} \,  \frac{\Gamma(i \, + \, 2 j \, + \alpha + 1)}{\Gamma(i + j + \alpha + 1)} 
\, G_{i}^{(\alpha - \beta + j \, , \, \beta + j)}(x)  \, , \nonumber  \\
\mbox{and } \ D^{j - 1}f_{N}(x) &= \ \sum_{i = -1}^{N - j} \frac{f_{i + j}}{|\| G_{i + j}^{(\beta \, , \, \alpha - \beta)} |\|^{2}}
\frac{\Gamma(i \, + \, 2 j \, +  \alpha)}{\Gamma(i + j + \alpha + 1)} G_{i + 1}^{(\beta + j - 1 \, , \, \alpha - \beta + j - 1)}(x) \, .
  \label{eqwert2}
\end{align}

Then, for $M > N$,
\begin{align}
&
\| D^{j}  \big( \frac{1}{\rho^{(\alpha - \beta \, , \, \beta)}(x)} 
\left( u_{M} \ - \ u_{N}  \right) \big) \|^{2}_{\rho^{(\alpha - \beta + j \, , \, \beta + j )}}
 \nonumber  \\
&= \ 
 \sum_{i = N - j + 1}^{M - j} c_{i + j}^{2} \, \left( \frac{\Gamma(i \, + \, 2 j \, + \alpha + 1)}{\Gamma(i + j + \alpha + 1)} \right)^{2}
\, |\| G_{i}^{(\alpha - \beta + j \, , \, \beta + j)} |\|^{2}   \nonumber \\
&= \ 
 \sum_{i = N - j + 1}^{M - j}  \frac{f_{i + j}^{2}}{ \lambda_{i + j}^{2} \, |\| G_{i + j}^{(\beta \, , \, \alpha - \beta)} |\|^{4}}
 \, \left( \frac{\Gamma(i \, + \, 2 j \, + \alpha + 1)}{\Gamma(i + j + \alpha + 1)} \right)^{2}
\, |\| G_{i}^{(\alpha - \beta + j \, , \, \beta + j)} |\|^{2}  \ \  \mbox{(using \eqref{defuNs2})} \nonumber \\
&\le \ 
C \,   \sum_{i = N - j + 1}^{M - j} 
\frac{f_{i + j}^{2}}{ |\| G_{i + j}^{(\beta \, , \, \alpha - \beta)} |\|^{4}} \, 
\left( \frac{\Gamma(i \, + \, 2 j \, + \alpha)}{\Gamma(i + j + \alpha + 1)} \right)^{2} \, 
|\| G_{i  + 1}^{(\beta + j - 1 \, , \, \alpha - \beta + j  - 1)} |\|^{2}   \ \ 
 \mbox{(using \eqref{errt2})} \nonumber \\
&= C \, \| D^{j - 1}  f_{M} \, - \, D^{j - 1} f_{N}  \|_{\rho^{(\beta + j - 1 \, , \, \alpha - \beta + j - 1 )}}^{2}  \ \ \ \ 
\mbox{(using \eqref{eqwert2})}  \, .  \label{thy1}
\end{align}

Assuming that $D^{j - 1}f \in L^{2}_{\rho^{(\beta + j - 1 \, , \, \alpha - \beta + j - 1)}}(\mrI)$, then $\{D^{j - 1} f_{n} \}$ is a Cauchy sequence in 
\linebreak[4]
$L^{2}_{\rho^{(\beta + j - 1 \, , \, \alpha - \beta + j - 1)}}(\mrI)$. Thus we can conclude that 
$D^{j}  \frac{1}{\rho^{(\alpha - \beta \, , \, \beta )}(x)} u(x)  \in  
L^{2}_{\rho^{(\alpha - \beta + j \, , \, \beta + j )}}(\mrI)$.  \\
\mbox{ } \hfill \qed
    
Combining Lemma \ref{thmrgqd} with the definition of $H^{l}_{\rho^{(a \, , \, b)} \,  , \, A}(\mrI)$, \eqref{defwHr}, we have the
following theorem.
\begin{theorem} \label{thmshift}
For $j \in \mathbb{N}$, if $f(x) \in H^{j-1}_{\rho^{(\beta \, , \, \alpha - \beta)} \,  , \, A}(\mrI)$,
then
$ \frac{1}{\rho^{(\alpha - \beta \, , \, \beta )}(x)} u(x) \,  \in H^{j}_{\rho^{(\alpha - \beta \, , \,  \beta)} \,  , \, A}(\mrI)$.
\end{theorem}
\mbox{ } \hfill \qed  \\

In the theory of linear differential equations a \textit{shift theorem} typically establishes that if the regularity of
the rhs function is improved by one order then the regularity of the solution also increases by one order.

 As $\rho^{(\alpha - \beta \, , \, \beta )}(x) \in C^{\infty}(a , b)$, for $0 < a < b < 1$, then 
 $\widehat{u}(x) \, := \, \frac{1}{\rho^{(\alpha - \beta \, , \, \beta )}(x)} u(x)$ will have the same regularity as
 $u(x)$ on $(a , b)$. Theorem \ref{thmshift} shows that 
 \underline{away from the endpoints} if the regularity of $f$ is improved by one order than the regularity of the
 solution also improves by one order.
 
 Also worthy of note is that, even though $\mbox{}_{RLC}\mcD_{r}^{\alpha}$ is a nonlocal operator, 
 Theorem \ref{thmshift} shows that $f$ may be singular at the endpoints without affecting the regularity
 of the solution away from the endpoints.

\underline{Summary of solution to \eqref{rlcde1}}: For $0 \le r \le 1$ and
$f \in H^{j-1}_{\rho^{(\beta \, , \, \alpha - \beta)} \,  , \, A}(\mrI)$ the RLC fractional
diffusion equation is well posed for all Dirichlet boundary conditions. The solution $\tilde{u}(x)$ is
decomposable into three pieces. Two pieces are explicitly determined by the values of the boundary
conditions, whereas the third piece $u(x)$ is determined by $f(x)$ and satisfies
$ \frac{1}{\rho^{(\alpha - \beta \, , \, \beta )}(x)} u(x) \,  \in H^{j}_{\rho^{(\alpha - \beta \, , \,  \beta)} \,  , \, A}(\mrI)$.

\subsection{Dirichlet and Neumann Boundary Condition}
\label{ssecDNBC}
Of interest in this section is the solution $\tilde{u}(x)$ of
\be
    \mbox{}_{RLC}\mcD_{r}^{\alpha} \tilde{u}(x) \ = \ f(x) \, , \ \ 0 < x < 1 \, , \ \ \ 
    \mbox{subject to } \mbox{}_{RLC}\mcF\tilde{u}(0) = A \, , \ \ \tilde{u}(1) = B \, .
\label{eqDNBC}
\ee

\underline{For $0 \le r < 1$}: \\
Consider $\tilde{u}(x) \ = \ u(x) \, + \, C_{1} k_{0}(x) \, + \, B$, with $u(x)$ given by Theorem \ref{thmuspg},
and $k_{0}(x)$ by \eqref{kRLC1}. Then $ \mbox{}_{RLC}\mcD_{r}^{\alpha} \tilde{u}(x) \ = \ f(x) \, , \ \ 0 < x < 1 \, , \ \  
    \mbox{and }  \ \tilde{u}(1) = B$.

Noting that, as $u(0) = u(1) = 0$, $\mbox{}_{RLC}\mcF u(x) \,  =  \,  \mbox{}_{RL}\mcF u(x)$, and using 
\eqref{NactG} and \eqref{chr111} we obtain
\be
 \mbox{}_{RLC}\mcF \tilde{u}(x) \ = \ - \, \sum_{i = 0}^{\infty} \mu_{i} \, c_{i} \, 
 G_{i+1}^{(\beta - 1 \, , \, \alpha - \beta - 1)}(x) \, + \, C_{1} \sigma_{0} \, .
\label{soldn}
\ee
Therefore,  $\mbox{}_{RLC}\mcF\tilde{u}((0) = A$ implies
\be
   C_{1} \ = \  \frac{1}{\sigma_{0}} \left( A \, + \,  \sum_{i = 0}^{\infty} \mu_{i} \, c_{i} \, 
 G_{i+1}^{(\beta - 1 \, , \, \alpha - \beta - 1)}(0) \right) \, , 
\label{defC1}
\ee
where convergence of the series is established in Lemma \ref{sumcvg1} for $0 \le r < 1$.

\underline{For $r = 1$}: (for which $\beta \, = \, \alpha - 1$) \\
For $f \in L^{2}_{\rho^{(\alpha - 1 \, , \, 1)}}(\mrI)$ we have that a necessary and sufficient
condition for $\tilde{u}(x)$ satisfying $ \mbox{}_{RLC}\mcD_{r}^{\alpha} \tilde{u}(x) \ = \ f(x) \, , \, x \in \mrI$, 
$\tilde{u}(1) \, = \, B$ to satisfy $\mbox{}_{RLC}\mcF\tilde{u}(0) = A$ is that the series in\eqref{defC1},
with $\beta \, = \, \alpha - 1$, converges. However, consider
\[
f(x) \ = \ \sum_{i = 2}^{\infty} \frac{f_{i}}{|\| G_{i}^{(\alpha - 1 \, , \, 1)} |\|^{2}}  G_{i}^{(\alpha - 1 \, , \, 1)}(x) \, , 
\ \ \mbox{where } \ f_{i} \ = \ (-1)^{i} \, \frac{1}{\log (i)} \, .
\]

Using \eqref{spm22g},
\[
|\| G_{i}^{(\alpha - 1 \, , \, 1)} |\|^{2} \ = \ \frac{1}{2 i \, + \, \alpha  + 1} 
\frac{\Gamma(i + \alpha) \, \Gamma(i + 2)}{\Gamma(i + 1) \, \Gamma(i + \alpha + 1)} \ = \ 
\frac{1}{2 i \, + \, \alpha  + 1} \, \frac{i + 1}{i + \alpha} \ \sim \ \frac{1}{2 i} \, .
\]

Note that
\begin{align*}
\| f \|_{L^{2}_{\rho^{(\alpha - 1 \, , \, 1)}}}^{2} \ = \ \sum_{i = 2}^{\infty} \frac{f_{i}^{2}}{|\| G_{i}^{(\alpha - 1 \, , \, 1)} |\|^{2}} 
&\sim \ 2 \, \int_{2}^{\infty} x \, \frac{1}{x^{2} \, ( \log(x) )^{2}} \, dx \\
&= \ 2 \, (-1) \left( \log(x) \right)^{-1} |_{x = 2}^{\infty} \ = \ \frac{2}{\log(2)} \ < \ \infty \, .
\end{align*}

However, corresponding to the series in  \eqref{defC1} we have
\begin{align*}
\sum_{i = 0}^{\infty} \mu_{i} \, c_{i} \, 
 G_{i+1}^{(\beta - 1 \, , \, \alpha - \beta - 1)}(0)  &= \
 \sum_{i = 2}^{\infty} \frac{(-1)^{i} \, f_{i}}{(i + \alpha) \, |\| G_{i}^{(\alpha - 1 \, , \, 1)} |\|^{2}} \\
&\sim \ \int_{2}^{\infty} \frac{1}{x \, \log(x) } \, dx \ = \ \log ( \log (x) ) |_{x = 2}^{\infty} \, \rightarrow \, \infty \, .
\end{align*}

Hence we conclude that for $r = 1$ and arbitrary $f \in L^{2}_{\rho^{(\beta , \alpha - \beta)}}(\mrI) \, = \,
L^{2}_{\rho^{(\alpha - 1 , 1)}}(\mrI)$ the problem \eqref{eqDNBC} is not well posed.

\underline{Summary of solution to \eqref{eqDNBC}}: 
For $0 \le r < 1$ and  $f \in H^{j-1}_{\rho^{(\beta \, , \, \alpha - \beta)} \,  , \, A}(\mrI)$ 
the RLC fractional
diffusion equation is well posed for mixed Dirichlet and Neumann boundary conditions. The solution $\tilde{u}(x)$ is
decomposable into three pieces. Two pieces are explicitly determined by the values of the boundary
conditions, whereas the third piece $u(x)$ is determined by $f(x)$ and satisfies
$ \frac{1}{\rho^{(\alpha - \beta \, , \, \beta )}(x)} u(x) \,  \in H^{j}_{\rho^{(\alpha - \beta \, , \,  \beta)} \,  , \, A}(\mrI)$. 
For $f \in H^{j-1}_{\rho^{(\alpha - 1 \, , \, 1)} \,  , \, A}(\mrI)$ and $r = 1$ problem \eqref{eqDNBC} is not well posed.

\subsection{Neumann Boundary Conditions}
\label{ssecNBC}
Of interest in this section is the solution $\tilde{u}(x)$ of
\be
    \mbox{}_{RLC}\mcD_{r}^{\alpha} \tilde{u}(x) \ = \ f(x) \, , \ \ 0 < x < 1 \, , \ \ \ 
    \mbox{subject to } \mbox{}_{RLC}\mcF\tilde{u}(0) = A \, , \ \ \mbox{}_{RLC}\mcF\tilde{u}(1) = B \, .
\label{eqNBC}
\ee

From integrating the differential equation we have
\be
\int_{0}^{1}  \mbox{}_{RLC}\mcD_{r}^{\alpha} \tilde{u}(s) \, ds \ = \ 
\mbox{}_{RLC}\mcF\tilde{u}(1) \, - \, \mbox{}_{RLC}\mcF\tilde{u}(0) \ = \ B \, - \, A \ = \ 
\int_{0}^{1} f(s) \, ds \, , 
\label{dnbc}
\ee
which gives the usual compatibility condition between the flux and the right hand side function
for a diffusion problem subject to Neumann boundary conditions.

\underline{For $0 < r < 1$}: \\
Assuming the compatibility condition is satisfied, from Section \ref{ssecDNBC} we have that,
for $C_{1}$ given by \eqref{defC1},
solutions to \eqref{eqNBC} are given by
\[
  \tilde{u}(x) \ = \ u(x) \, + \, C_{1}  k_{0}(x) \, + \, C_{3} \, ,
\]
for any $C_{3} \in \real$.

From  \eqref{soldn},
\begin{align*}
\mbox{}_{RLC}\mcF\tilde{u}(1) \, - \, \mbox{}_{RLC}\mcF\tilde{u}(0) &= \ B \, - \, A \ = \
- \sum_{i = 0}^{\infty} \mu_{i} \, c_{i} \, \left( G_{i+1}^{(\beta - 1 \, , \, \alpha - \beta - 1)}(1) \, - \,
   G_{i+1}^{(\beta - 1 \, , \, \alpha - \beta - 1)}(0) \right)  \\
&= \ - \sum_{i = 0}^{\infty} \frac{\mu_{i}}{\lambda_{i}} \frac{f_{i}}{|\| G_{i}^{(\beta , \alpha - \beta)} |\|^{2}} \,   
\int_{0}^{1} \frac{d}{ds}  G_{i+1}^{(\beta - 1 \, , \, \alpha - \beta - 1)}(s) \, ds  \\
&= \ \int_{0}^{1}  \sum_{i = 0}^{\infty} \frac{\mu_{i}}{\lambda_{i}} \frac{f_{i}}{|\| G_{i}^{(\beta , \alpha - \beta)} |\|^{2}} \,   
\frac{\Gamma(i + 1 + \alpha)}{\Gamma(i + \alpha)} \,  G_{i}^{(\beta  \, , \, \alpha - \beta)}(s) \, ds  \\
&= \ \int_{0}^{1} f(s) \, ds \, ,
\end{align*}
confirming the compatibility condition.

\underline{For $r = 0$ and $r = 1$}: \\
Analogous to the discussion for the mixed boundary condition problem discussed in Section \ref{ssecDNBC}, \\
For $r = 0$: (for which $\beta = 1$) For $f \in L^{2}_{\rho^{(1 , \alpha - 1)}}(\mrI)$ the problem \eqref{eqNBC} is not well
posed.  \\
For $r = 1$: (for which $\beta \, = \, \alpha - 1$) For $f \in L^{2}_{\rho^{(\alpha - 1 \, , \, 1)}}(\mrI)$ 
the problem \eqref{eqNBC} is not well
posed.  \\

\underline{Summary of solution to \eqref{eqNBC}}: 
For $0 < r < 1$ and $f \in H^{j-1}_{\rho^{(\beta \, , \, \alpha - \beta)} \,  , \, A}(\mrI)$ the RLC fractional
diffusion equation is well posed for  Neumann boundary conditions, subject to the boundary
conditions satisfying the usual compatibility condition \eqref{dnbc}. The solution is only determined up to an
additive solution. Additionally,
the solution $\tilde{u}(x)$ is
decomposable into three pieces, an undetermined constant, 
a piece explicitly determined by the values of the boundary conditions, 
and a third piece $u(x)$ determined by $f(x)$ satisfying
$ \frac{1}{\rho^{(\alpha - \beta \, , \, \beta )}(x)} u(x) \,  \in H^{j}_{\rho^{(\alpha - \beta \, , \,  \beta)} \,  , \, A}(\mrI)$. 
For $r = 0$ and $f \in L^{2}_{\rho^{(1 , \alpha - 1)}}(\mrI)$, or $r = 1$ and $f \in L^{2}_{\rho^{(\alpha - 1 \, , \, 1)}}(\mrI)$ 
the problem \eqref{eqNBC} is not well posed.

\setcounter{equation}{0}
\setcounter{figure}{0}
\setcounter{table}{0}
\setcounter{theorem}{0}
\setcounter{lemma}{0}
\setcounter{corollary}{0}
\section{Existence and Regularity of the RL Fractional Diffusion Model}
\label{sec_exRL}
In this section we investigate the existence and regularity of the steady state RL fractional diffusion model subject to
various boundary conditions.

From Lemma \ref{ker1r} we have that
\begin{align}
 ker( \mbox{}_{RL}\mcD_{r}^{\alpha} ) &= \ span \left\{  (1 - x)^{\alpha - \beta - 1} \, x^{\beta - 1} \ , \   
 (1 - x)^{\alpha - \beta - 1} \, x^{\beta}  \right\}  \nonumber \\
&= \ span \{ (1 - x)^{\alpha - \beta - 1} \, x^{\beta}  \ , 
\ (1 - x)^{\alpha - \beta} \, x^{\beta - 1} \}
     \label{kRL1}
\end{align}

The singular endpoint behavior of the kernel at both endpoints is more apparent in representation \eqref{kRL1}.

From \cite{erv162} we have, as $\rho^{(\alpha - \beta , \beta)}(0) \, = \, \rho^{(\alpha - \beta , \beta)}(1) \, = \, 0$,
\[
 \mbox{}_{RL}\mcD_{r}^{\alpha}  \rho^{(\alpha - \beta , \beta)}(x) G_{n}^{(\alpha - \beta , \beta)}(x) \ = \ 
 \mbox{}_{RLC}\mcD_{r}^{\alpha} \rho^{(\alpha - \beta , \beta)}(x) G_{n}^{(\alpha - \beta , \beta)}(x) 
 \ = \ \lambda_{n} G_{n}^{(\beta , \alpha - \beta)}(x)  \, .
\]

For  $ \mbox{}_{RL}\mcD_{r}^{\alpha} \cdot $ we also have the following.
\begin{theorem} \label{thmLLJ} 
Under  \textbf{Condition A}, 
for $n=0,1,2,\cdots,$
\begin{align}
    \mbox{}_{RL}\mcD_{r}^{\alpha}  \rho^{(\alpha - \beta - 1 , \beta - 1)}(x) G_{n}^{(\alpha - \beta - 1 , \beta - 1)}(x)
    &= \  \kappa_{n} \,  G_{n - 2}^{(\beta + 1 , \alpha - \beta + 1)}(x)
    \, ,   \label{lljt1}  \\
\mbox{ where  } \ \kappa_{n} \ = \  
c_{*}^{*} \,  \frac{\Gamma(n + \alpha + 1)}{\Gamma(n + 1)} \,\quad 
&\mbox{ and } G_{j}^{(\cdot , \cdot)}(x) \, = \, 0 \, , \ \ \mbox{ for } j < 0  \, . \nonumber 
\end{align}
\end{theorem}
\textbf{Proof}: 
From \eqref{eqC4} we have
\begin{equation}  \label{jacobi_deri}
\frac{d^2}{dx^2}G_n^{(\beta-1 , \alpha-\beta-1)}(x) \ = \ \frac{\Gamma(n + \alpha +1)}{\Gamma(n + \alpha - 1)}
G_{n-2}^{(\beta+1 , \alpha-\beta+1)}(x) \, .
\end{equation}
Combining \eqref{jacobi_deri} and Lemma \ref{ch_lma} we obtain \eqref{lljt1}. \\
\mbox{ } \hfill \qed

Note that for $n = 0$ and $1$, $ \rho^{(\alpha - \beta - 1 , \beta - 1)}(x) G_{n}^{(\alpha - \beta - 1 , \beta - 1)}(x) 
\in ker( \mbox{}_{RL}\mcD_{r}^{\alpha} )$.

For $f(x) \in L^{2}_{ \rho^{(\beta , \alpha - \beta)}}(\mrI)$, $u(x) \in L^{2}_{ \rho^{(-(\alpha -\beta) , -\beta)}}(\mrI)$
satisfying 
 $\mbox{}_{RL}\mcD_{r}^{\alpha}  u(x) \ = \ f(x)$ can be expressed as given in Theorem \ref{thmuspg}.
 
 Using Theorem \ref{thmLLJ}, in a similar fashion as was done in \cite{erv162}, 
 $f(x) \in L^{2}_{ \rho^{(\beta + 1, \alpha - \beta + 1)}}(\mrI)$, $w(x) \in L^{2}_{ \rho^{(-(\alpha -\beta - 1) , -(\beta - 1) )}}(\mrI)$
satisfying 
 $\mbox{}_{RL}\mcD_{r}^{\alpha}  w(x) \ = \ f(x)$ can be expressed as
\begin{align}
w(x) &= \  \rho^{( \alpha -\beta - 1  , \beta - 1 )}(x) \, \sum_{i = 2}^{\infty} w_{i} \, G_{i}^{( \alpha -\beta - 1  , \beta - 1 )}(x)
\, , \ \ \ \ \mbox{  where  }   \label{repHRL}  \\
 w_{i} &= \  \frac{1}{\kappa_{i}}  \, 
 \frac{1}{|\| G_{i}^{(\beta + 1, \alpha - \beta + 1)} |\|^{2}} \, 
 \int_{0}^{1} \rho^{(\beta + 1, \alpha - \beta + 1)}(s) \, G_{i}^{(\beta + 1, \alpha - \beta + 1)}(s) \, f(s) \, ds \, .
 \nonumber 
\end{align}

Note that  $L^{2}_{ \rho^{(\beta, \alpha - \beta )}}(\mrI) \subset L^{2}_{ \rho^{(\beta  + 1, \alpha - \beta + 1)}}(\mrI)$.

In order to contrast the solutions of the $RLC$ and $RC$ diffusion equations, in this section we will assume
that $f(x) \in L^{2}_{ \rho^{(\beta , \alpha - \beta)}}(\mrI)$.

\subsection{Dirichlet Boundary Conditions}
\label{ssec_RLdbc}
For $f(x) \in L^{2}_{ \rho^{(\beta , \alpha - \beta)}}(\mrI)$, we consider the RL diffusion equation with Dirichlet boundary
conditions: 
\be
 \mbox{}_{RL}\mcD_{r}^{\alpha}  \tilde{u}(x) \ = \ f(x) \, \ \ 0 < x < 1 \, , \ \ \mbox{ subject to }
  \tilde{u}(0) \, = \, A \, , \ \ \tilde{u}(1) \, = \, B \, .
  \label{eqtu1}
\ee
Using \eqref{kRL1} and Theorem \ref{thmuspg} the general solution of \eqref{eqtu1} can be expressed as
\[
 \tilde{u}(x) \ = \ C_{1} \, (1 - x)^{\alpha - \beta - 1} x^{\beta} \ + \ C_{2}  \, (1 - x)^{\alpha - \beta} x^{\beta - 1} 
 \ + \ \rho^{(\alpha - \beta , \beta)}(x) \, \sum_{j = 0}^{\infty} c_{j} \, G_{j}^{(\alpha - \beta , \beta)}(x) \, .
\]
Now,
\be
\tilde{u}(0) \, = \, A \ \Longrightarrow \ A \, = \, C_{2} \,  \lim_{x \rightarrow 0} x^{\beta - 1} \, , \ \ 
\mbox{ and } \ \  
\tilde{u}(1) \, = \, B \ \Longrightarrow \ B \, = \, C_{1} \, \lim_{x \rightarrow 1} (1 - x)^{\alpha - \beta - 1} \, .
 \label{eqtu2}
\ee
For $A , \, B \in \real$, in order that \eqref{eqtu2} defines a finite value for $C_{1}$ and $C_{2}$ we must have
$A \, = \, B \, = \, 0 \ \Longrightarrow \ C_{1} \, = \, C_{2} \, = \, 0$. Recall that in the case of homogeneous
Dirichlet boundary conditions problems \eqref{DefDRLC} and \eqref{DefDRL}  coincide.

In place of \eqref{eqtu1}, if we consider the problem:
\be
 \mbox{}_{RL}\mcD_{r}^{\alpha}  \tilde{u}(x) \ = \ f(x) \, \ \ 0 < x < 1 \, , \ \ \mbox{ subject to }
 \lim_{x \rightarrow 0} \tilde{u}(x) \, = \, A x^{\beta - 1}\, , \ \  \lim_{x \rightarrow 1} \tilde{u}(x) \, = \, B  (1 - x)^{\alpha - \beta - 1} \, ,
  \label{eqtu3}
\ee
then the solution is well defined, satisfying
\be
\tilde{u}(x) \ = \ A (1 - x)^{\alpha - \beta} \, x^{\beta - 1} \, + \, B  (1 - x)^{\alpha - \beta - 1} x^{\beta} \, + \, 
\rho^{( \alpha - \beta , \beta )}(x) \, \sum_{i = 0}^{\infty} c_{i} \, G_{i}^{( \alpha - \beta , \beta )}(x) \, .
\label{eqtu4}
\ee

\underline{Summary of solution to \eqref{eqtu1}}: 
For $0 \le r \le 1$ and $f \in H^{j-1}_{\rho^{(\beta \, , \, \alpha - \beta)} \,  , \, A}(\mrI)$ in order for the
the RL fractional diffusion equation to be well posed the solution must have a specific, prescribed singular
behavior at the endpoints of the interval. In that case, the solution $\tilde{u}(x)$ is
decomposable into three pieces, two singular pieces that are determined by the boundary conditions and a regular 
piece determined by $f(x)$. This regular piece is the same as discussed in Theorem \ref{thmshift}.
 
\subsection{Dirichlet and Neumann Boundary Condition}
\label{ssecDNBCRL}
In this section we consider the problem:
\be
    \mbox{}_{RL}\mcD_{r}^{\alpha} \tilde{u}(x) \ = \ f(x) \, , \ \ 0 < x < 1 \, , \ \ \ 
    \mbox{subject to } \mbox{}_{RL}\mcF\tilde{u}(0) = A \, , 
    \ \ \lim_{x \rightarrow 1} \tilde{u}(x) \, = \, B  (1 - x)^{\alpha - \beta - 1} \, .
\label{eqtu5}
\ee
\underline{For $0 \le r < 1$}: \\
It is convenient to express the solution as
\be
\tilde{u}(x) \ = \ 
\tilde{u}(x) \ = \ C_{1} (1 - x)^{\alpha - \beta - 1} \, x^{\beta} \, + \, C_{3}  (1 - x)^{\alpha - \beta - 1} x^{\beta - 1} \, + \, 
\rho^{( \alpha - \beta , \beta )}(x) \, \sum_{i = 0}^{\infty} c_{i} \, G_{i}^{( \alpha - \beta , \beta )}(x) \, .
\label{poif1}
\ee
Using Lemma \ref{ker1r} and Lemma \ref{lmaegp},
\[
\mbox{}_{RL}\mcF\tilde{u}(x) \ = \ 
- C_{1} \mu_{-1} \ + \ 0 \ - \ \sum_{i = 0}^{\infty} \mu_{i} \, c_{i} \, G_{i}^{(\beta - 1,  \alpha - \beta - 1 )}(x) \, .
\]
Therefore $\mbox{}_{RL}\mcF\tilde{u}(0) = A$ implies
\be
C_{1} \ = \ - \frac{1}{\mu_{-1}} \left( A \, + \, 
  \sum_{i = 0}^{\infty} \mu_{i} \, c_{i} \, G_{i}^{(\beta - 1,  \alpha - \beta - 1 )}(0)  \right) \, ,
\label{eqtu6}
\ee
where convergence of the series is established in Lemma \ref{sumcvg1} for $0 \le r < 1$.

The boundary condition at $x = 1$ implies
\[ 
(C_{1} \, + \, C_{3} ) \, \lim_{x \rightarrow 1}   (1 - x)^{\alpha - \beta - 1} 
\ = \ B \, \lim_{x \rightarrow 1} (1 - x)^{\alpha - \beta - 1}  \, ,
\]
which gives $C_{3} \ = \ B \, - \, C_{1}$. The solution is then given by
\[
\tilde{u}(x) \ = \ C_{1} (1 - x)^{\alpha - \beta - 1} \, x^{\beta} \, + \, 
(B  -  C_{1})  (1 - x)^{\alpha - \beta - 1} x^{\beta - 1} \, + \, 
\rho^{( \alpha - \beta , \beta )}(x) \, \sum_{i = 0}^{\infty} c_{i} \, G_{i}^{( \alpha - \beta , \beta )}(x) \, ,
\]
with $C_{1}$ given by \eqref{eqtu6}.

\underline{For $r = 1$}: (for which $\beta \, = \, \alpha - 1$) \\
Analogous to the discussion in Section \ref{ssecDNBC}, \eqref{eqtu5} is not well posed for $r = 1$ and 
arbitrary $f \in L^{2}_{\rho^{(\alpha - 1 \, , \, 1)}}(\mrI)$.

\underline{Summary of solution to \eqref{eqtu5}}: For $0 \le r < 1$,
$f \in H^{j-1}_{\rho^{(\beta \, , \, \alpha - \beta)} \,  , \, A}(\mrI)$,
a flux boundary condition imposed at $x = 0$, and a prescribed boundary condition behavior of the
form $(1 - x)^{\alpha - \beta - 1}$ at $x = 1$, the RL fractional
diffusion equation \eqref{eqtu5} is well posed. For the case $r = 1$, and 
arbitrary $f \in L^{2}_{\rho^{(\alpha - 1 \, , \, 1)}}(\mrI)$, \eqref{eqtu5} is not well posed.

\subsection{Neumann Boundary Conditions}
\label{ssecRLNBC}
Of interest in this section is the solution $\tilde{u}(x)$ of
\be
    \mbox{}_{RL}\mcD_{r}^{\alpha} \tilde{u}(x) \ = \ f(x) \, , \ \ 0 < x < 1 \, , \ \ \ 
    \mbox{subject to } \mbox{}_{RLC}\mcF\tilde{u}(0) = A \, , \ \ \mbox{}_{RLC}\mcF\tilde{u}(1) = B \, .
\label{eqNRL}
\ee

Again we have the usual compatibility condition between the flux and the right hand side function
\be
\int_{0}^{1}  \mbox{}_{RL}\mcD_{r}^{\alpha} \tilde{u}(s) \, ds \ = \ 
\mbox{}_{RL}\mcF\tilde{u}(1) \, - \, \mbox{}_{RL}\mcF\tilde{u}(0) \ = \ B \, - \, A \ = \ 
\int_{0}^{1} f(s) \, ds \, ,
\label{dnbc2}
\ee
which we assume is satisfied.

\underline{For $0 < r < 1$}: \\
From Section \ref{ssecDNBCRL} we have that the
solution to \eqref{eqNRL} is given by \eqref{poif1},  with $C_{1}$ determined by \eqref{eqtu6}
and $C_{3} \in \real$ an arbitrary constant. 

\underline{For $r = 0$ and $r = 1$}: \\
Analogous to the discussion for the  $\mbox{}_{RLC}\mcD_{r}^{\alpha}$ operator in
discussed in Section \ref{ssecNBC}, \\
For $r = 0$: (for which $\beta = 1$) For $f \in L^{2}_{\rho^{(1 , \alpha - 1)}}(\mrI)$ the problem \eqref{eqNRL} is not well
posed.  \\
For $r = 1$: (for which $\beta \, = \, \alpha - 1$) For $f \in L^{2}_{\rho^{(\alpha - 1 \, , \, 1)}}(\mrI)$ 
the problem \eqref{eqNRL} is not well
posed.  \\

\underline{Summary of solution to \eqref{eqNRL}}: 
For $0 < r < 1$ and $f \in H^{j-1}_{\rho^{(\beta \, , \, \alpha - \beta)} \,  , \, A}(\mrI)$ 
the RL fractional
diffusion equation is well posed for  Neumann boundary conditions, subject to the boundary
conditions satisfying the usual compatibility condition \eqref{dnbc2}. The solution is only determined up to an
additive solution. Additionally, the solution $\tilde{u}(x)$ is
decomposable into three pieces, an undetermined constant, 
a piece explicitly determined by the values of the boundary conditions, 
and a third piece $u(x)$ determined by $f(x)$ satisfying
$ \frac{1}{\rho^{(\alpha - \beta \, , \, \beta )}(x)} u(x) \,  \in H^{j}_{\rho^{(\alpha - \beta \, , \,  \beta)} \,  , \, A}(\mrI)$. 
For $r = 0$ and $f \in L^{2}_{\rho^{(1 , \alpha - 1)}}(\mrI)$, or $r = 1$ and $f \in L^{2}_{\rho^{(\alpha - 1 \, , \, 1)}}(\mrI)$ 
the problem \eqref{eqNRL} is not well posed.

\setcounter{equation}{0}
\setcounter{figure}{0}
\setcounter{table}{0}
\setcounter{theorem}{0}
\setcounter{lemma}{0}
\setcounter{corollary}{0}
\section{Conclusions}
In this paper we have investigated the well posedness and regularity of the solution to fractional diffusion equations
\eqref{DefDRLC} and \eqref{DefDRL}. In the case of homogeneous Dirichlet boundary conditions or Neumann
boundary conditions the solutions to \eqref{DefDRLC} and \eqref{DefDRL} agree. However, for nonhomogeneous
Dirichlet boundary conditions that is not the case. Specifically, for nonhomogeneous
Dirichlet boundary conditions the solution to \eqref{DefDRLC} is bounded on $(0 , 1)$, whereas for 
the problem \eqref{DefDRL} to be well posed specific singular behavior at the endpoints must be specified.
Regarding the regularity of the solution, we have shown that the solution, away from the endpoints, satisfies
a shift theorem with respect to the regularity of the rhs function.

\appendix

\setcounter{equation}{0}
\setcounter{figure}{0}
\setcounter{table}{0}
\setcounter{theorem}{0}
\setcounter{lemma}{0}
\setcounter{corollary}{0}
\section{Ancillary Properties and Proofs}
\label{apdxAPP}
In this section we presents some ancillary results used in establishing the existence and regularity properties given above.

\begin{lemma} \label{spec_poly}
Under \textbf{Condition A}, 
for $n=0,1,2,\cdots,$
\begin{equation}\label{spec_poly_r}
\begin{array}{rl}
\textbf{I}_r^{2-\alpha}(1-x)^{\alpha-\beta-1}x^{\beta-1}x^n=&\sum\limits_{k=0}^na_{n,k}x^k,
\end{array}
\end{equation}
where $a_{n,k}=(-1)^{n+1} \, c_{*}^{*} \, \Gamma(\alpha-\beta)
\frac{(-1)^k \, \Gamma(\alpha-1+k)}{\Gamma(\alpha-\beta-n+k) \, \Gamma(n+1-k) \, \Gamma(k+1)}$.
\end{lemma}
\textbf{Proof}: 
With $u(x)=(1-x)^{\alpha-\beta-1}x^{\beta-1}x^n,$ using Maple we obtain that
\begin{equation}
\begin{array}{rl}
D^{-(2-\alpha)}u(x)&=\frac{\Gamma(\beta+n)}{\Gamma(2-\alpha+\beta+n)}x^{n+1-\alpha+\beta}
{_2F_1}(n+\beta,\beta-\alpha+1;2-\alpha+\beta+n,x) \, ,
\end{array}
\end{equation}
and
\begin{equation}
\begin{array}{rl}
D^{-(2-\alpha)*}u(x)&=\frac{\Gamma(\alpha-\beta-n-1)}{\Gamma(1-\beta-n)}x^{n+1-\alpha+\beta}{_2F_1}(n+\beta,\beta-\alpha+1;2-\alpha+\beta+n,x)\\
&+(-1)^{n+1}\Gamma(\alpha-\beta)\sum\limits_{k=0}^{n}\frac{(-1)^k\csc(\pi(\alpha-\beta)+k\pi)\sin(\pi\alpha+k\pi)\Gamma(\alpha-1+k)}
{\Gamma(\alpha-\beta-n+k)\Gamma(n+1-k)\Gamma(k+1)}x^k \, ,
\end{array}
\end{equation}
where ${_2F_1}(\cdot , \cdot ; \cdot , x)$ denotes the Gaussian three parameter hypergeometric function.

Using the identity
\[
\Gamma(1-z)=\frac{\pi}{\sin(\pi z)}\frac{1}{\Gamma(z)},
\]
it follows that
\begin{equation}\label{iden_1}
\Gamma(1-\beta-n)=\frac{(-1)^n\pi}{\sin(\pi\beta)}\frac{1}{\Gamma(\beta+n)} \, ,
\end{equation}
\begin{equation}\label{iden_2}
\mbox{and } \ \ \Gamma(2-\alpha+\beta+n)=\frac{(-1)^{n+1}\pi}{\sin(\pi(\alpha-\beta))}\frac{1}{\Gamma(\alpha-\beta-n-1)} \, .
\end{equation}
From (\ref{iden_1}) and (\ref{iden_2}) we obtain,
\begin{equation} \label{iden_12}
\frac{\Gamma(\beta+n)}{\Gamma(2-\alpha+\beta+n)}=-\frac{\sin(\pi(\alpha-\beta))}{\sin(\pi\beta)}
\frac{\Gamma(\alpha-\beta-n-1)}{\Gamma(1-\beta-n)}.
\end{equation}

Using \eqref{iden_12}, the coefficient of $x^{n+1-\alpha+\beta}{_2F_1(\cdot)}$ in the linear combination  \linebreak[4]
$(rD^{-(2-\alpha)}+(1-r)D^{-(2-\alpha)*})u(x)$ is
\begin{align*}
& r \, \frac{\Gamma(\beta+n)}{\Gamma(2-\alpha+\beta+n)} \, + \, (1 - r)\frac{\Gamma(\alpha-\beta-n-1)}{\Gamma(1-\beta-n)} 
\ = \ \frac{\Gamma(\alpha-\beta-n-1)}{\Gamma(1-\beta-n)} \left( - r\frac{\sin(\pi(\alpha-\beta))}{\sin(\pi\beta)} \, + \, 1-r \right) \\
& \quad = \ \frac{\Gamma(\alpha-\beta-n-1)}{\Gamma(1-\beta-n)}
 \left( -  \, \frac{\sin(\pi\beta)}{\sin(\pi(\alpha-\beta))+\sin(\pi\beta)} \frac{\sin(\pi(\alpha-\beta))}{\sin(\pi\beta)} 
 \, + \, \frac{\sin(\pi(\alpha-\beta))}{\sin(\pi(\alpha-\beta))+\sin(\pi\beta)} \right) \\
& \quad = \ 0 \, .
\end{align*}

Next, it is straightforward to show that
\begin{equation*}
\csc(\pi(\alpha-\beta)+k\pi)\sin(\pi\alpha+k\pi)=\frac{\sin(\pi\alpha)}{\sin(\pi(\alpha-\beta))},
\end{equation*}
Then, as $(1 - r) \, \sin(\pi \alpha) / \sin(\pi (\alpha - \beta)) \ = \ c_{*}^{*}$, we obtain,
\begin{equation*}
\textbf{I}_r^{2-\alpha}u(x) \ = \ (-1)^{n+1} \, c_{*}^{*} \, \Gamma(\alpha-\beta) \, \sum\limits_{k=0}^n
\frac{(-1)^k \, \Gamma(\alpha-1+k)}{\Gamma(\alpha-\beta-n+k) \, \Gamma(n+1-k) \, \Gamma(k+1)} x^k \, .
\end{equation*}
\mbox{ } \hfill \qed

\textbf{Remark}: Note that in an analogous manner as used in Lemma \ref{spec_poly}, we have that
\begin{equation}  \label{erwt19}
\textbf{I}_{1-r}^{2-\alpha} \, (1-x)^{\beta-1}x^{\alpha-\beta-1}x^n \ = \ \sum\limits_{k=0}^{n}b_{n,k}x^k,
\end{equation}
where $b_{n,k}=(-1)^{n+1} \, c_{*}^{*} \, 
\Gamma(\beta)\frac{(-1)^k\Gamma(\alpha-1+k)}{\Gamma(\beta-n+k)\Gamma(n+1-k)\Gamma(k+1)}.$

\begin{lemma}\label{ch_lma}
Under \textbf{Condition A}, 
for $n=0,1,2,\cdots,$
\begin{equation} \label{spec_integral}
\textbf{I}_r^{2-\alpha}\rho^{(\alpha-\beta-1,\beta-1)}(x) \, G_n^{(\alpha-\beta-1,\beta-1)}(x) \ = \
\sigma_n \, G_n^{(\beta-1,\alpha-\beta-1)}(x),
\end{equation}
where $\sigma_{n}$ is given by \eqref{chr113}.
\end{lemma}
\textbf{Proof}: 
Using the orthogonality of $\left\{ G_n^{(\alpha-\beta-1,\beta-1)}(x) \right\}_{n = 0}^{\infty}$ 
with respect to the weight function 
$\rho^{(\alpha-\beta-1,\beta-1)}(x)$, we have
that for any $p(x)\in \mathcal{P}_{n-1}(x)$, 
\be
\left( G_n^{(\alpha-\beta-1,\beta-1)}(x) \, , \, p(x) \right)_{\rho^{(\alpha-\beta-1,\beta-1)}} \ = \ 0 \, .
\label{eqgd1}
\ee
Additionally, up to a constant, \eqref{eqgd1}, defined the $n^{th}$ order polynomial $G_n^{(\alpha-\beta-1,\beta-1)}(x)$.

For $p(x)\in \mathcal{P}_{n-1}(x),$ from \eqref{erwt19}, there exist $\hat{p}(x)\in \mathcal{P}_{n-1}(x)$ such that 
\be
\textbf{I}_{1-r}^{2-\alpha}\rho^{(\beta-1,\alpha-\beta-1)}(x)p(x)=\hat{p}(x).
\label{eqgd2}
\ee
Then
\begin{align*}
& \left(\textbf{I}_r^{2-\alpha} \rho^{(\alpha-\beta-1,\beta-1)}(x) \, G_n^{(\alpha-\beta-1,\beta-1)}(x) \, , \, 
p(x) \right)_{\rho^{(\beta-1,\alpha-\beta-1)}}   \nonumber \\
& \quad \quad = \ \left(\textbf{I}_r^{2-\alpha} \rho^{(\alpha-\beta-1,\beta-1)}(x) \, G_n^{(\alpha-\beta-1,\beta-1)}(x) \, , \, 
 \rho^{(\beta-1,\alpha-\beta-1)} \, p(x) \right)  \nonumber  \\
& \quad \quad = \ \left(\rho^{(\alpha-\beta-1,\beta-1)}(x) \, G_n^{(\alpha-\beta-1,\beta-1)}(x) \, , \, 
 \textbf{I}_{1-r}^{2-\alpha} \rho^{(\beta-1,\alpha-\beta-1)} \, p(x) \right)  \nonumber   \\
& \quad \quad = \ \left(\rho^{(\alpha-\beta-1,\beta-1)}(x) \, G_n^{(\alpha-\beta-1,\beta-1)}(x) \, , \, 
 \hat{p}(x) \right)  
 \ \ \ \mbox{(using \eqref{eqgd2})}   \nonumber   \\ 
& \quad \quad = \ 0 \, ,  \ \ \ \mbox{(using \eqref{eqgd1})}. 
\end{align*}
Hence $\textbf{I}_r^{2-\alpha}\rho^{(\alpha-\beta-1,\beta-1)}(x)G_n^{(\alpha-\beta-1,\beta-1)}(x)=CG_n^{(\beta-1,\alpha-\beta-1)}(x)$ for $C\in \real$.\\
As the coefficient of $x^n$ in $G_n^{(\alpha-\beta-1,\beta-1)}(x)$ and $G_n^{(\beta-1,\alpha-\beta-1)}(x)$ is the same, 
then from Lemma \ref{spec_poly},
\[
C \ = \ - c_{*}^{*} \, \frac{\Gamma(n + \alpha - 1)}{\Gamma(n + 1)} \ = \sigma_n \, .
\]
\mbox{ } \hfill \qed

The bound obtained in the following lemma is used in establishing the regularity of $Du$ in
Corollary \ref{corRegDu}  in Section \ref{ssec_RLCdbc}.
\begin{lemma} \label{lmageq}
For $j \, = \, 0, 1, 2, \ldots$
\begin{equation}
 \frac{1}{2} \, \le \, \frac{ |\| G_{j}^{(\alpha - \beta \, , \, \beta)} |\|^{2} }{ |\| G_{j + 1}^{(\beta - 1 \, , \, \alpha - \beta - 1)} |\|^{2} } 
  \, = \, \frac{j + 1}{j + \alpha} \, \le \, 1 \, .
 \label{gge2}
\end{equation}
\end{lemma}
\textbf{Proof}: From \eqref{spm22g},
\begin{align}
\frac{ |\| G_{j}^{(\alpha - \beta \, , \, \beta)} |\|^{2} }{ |\| G_{j + 1}^{(\beta - 1 \, , \, \alpha - \beta - 1)} |\|^{2} } 
&= \ \frac{1}{2 j \, + \, \alpha \, + \, 1} \, 
\frac{\Gamma(j + \alpha - \beta + 1) \, \Gamma(j + \beta + 1)}{\Gamma(j + 1) \, \Gamma(j + \alpha + 1)} \ \
\frac{2 j \, + \, \alpha \, + \, 1}{1} \, 
\frac{\Gamma(j + 2) \, \Gamma(j + \alpha + 1)}{\Gamma(j + \beta + 1) \, \Gamma(j + \alpha - \beta + 1)}   \nonumber \\
&= \ \frac{j + 1}{j + \alpha} \le 1 \, .      \label{nrto1}
\end{align}
\mbox{ } \hfill \qed


The following lemma is used in the proof of Lemma \ref{thmrgqd} in Section \ref{ssec_RLCdbc}.
\begin{lemma} \label{lmahf1}
For $j \in \mathbb{N}$, there exists $C > 0$ such that
\begin{equation}
 \frac{(i + 2j + \alpha)^{2}}{\lambda_{i + j}^{2}} \,
 \frac{ |\| G_{i}^{(\alpha - \beta + j \, , \, \beta + j)} |\|^{2}}{ |\| G_{i + 1}^{(\beta + j - 1 \, , \, \alpha -  \beta + j - 1)} |\|^{2}}
 \ \le \ C \, .
 \label{errt2}
\end{equation}
\end{lemma}
\textbf{Proof}: From \eqref{nmeqG} and \eqref{spm22g} ,
\begin{align}
 \frac{ |\| G_{i}^{(\alpha - \beta + j \, , \, \beta + j)} |\|^{2}}{ |\| G_{i + 1}^{(\beta + j - 1 \, , \, \alpha -  \beta + j - 1)} |\|^{2}}
 &= \ 
\frac{ |\| G_{i}^{(\alpha - \beta - j \, , \, \beta - j)} |\|^{2}}{ |\| G_{i + 1}^{( \alpha - \beta + j - 1 \, , \, \beta + j - 1)} |\|^{2}}
 \nonumber \\
 & = \  \frac{1}{(2 i \, + \, \alpha \, + \, 2 j \, + \, 1)} \, 
  \frac{\Gamma( i + \alpha - \beta + j + 1) \, \Gamma(i + \beta + j + 1)}{\Gamma(i + 1) \, \Gamma(i + \alpha \,  + \, 2 j \, + 1)}  
  \nonumber \\
& \quad \quad \quad 
  \cdot \, (2 i \, + \, \alpha \, + \, 2 j \, + \, 1) \, 
  \frac{\Gamma( i  + 2) \, \Gamma(i + \alpha \, + \,  2 j)}{\Gamma(i + \alpha - \beta  +  j + 1) \, \Gamma(i + \beta + j + 1) }  
  \nonumber \\
&= \ 
  \frac{(i + 1)}{(i + \alpha \,  + \, 2 j)}  . \label{ppol1}
\end{align}

Using Stirling's formula,
\begin{align}
 \frac{1}{| \lambda_{i + j} |} \ = \ C \ \frac{\Gamma(i + j + 1)}{\Gamma(i + j + \alpha + 1)}
  &\sim 
   \ \left( i + j + 1 \right)^{- \alpha} \ \sim \ i^{- \alpha}  \, . \label{uy3}   \end{align}
Combining \eqref{ppol1} and \eqref{uy3} we obtain
\begin{align*}
  \frac{(i + 2j + \alpha)^{2}}{\lambda_{i + j}^{2}} \,
 \frac{ |\| G_{i}^{(\alpha - \beta + j \, , \, \beta + j)} |\|^{2}}{ |\| G_{i + 1}^{(\beta + j - 1 \, , \, \alpha -  \beta + j - 1)} |\|^{2}}  
  &\sim 
i^{- 2 \alpha} \, ( i + 2j + \alpha)^{2} \, \frac{(i + 1)}{(i + \alpha \,  + \, 2 j)}  \\
&\sim \ i^{- 2 (\alpha - 1)} \ \rightarrow \, 0 \, , \quad \mbox{ as } i \rightarrow \infty \, ,
\end{align*}
from which \eqref{errt2} follows.  \\
\mbox{ } \hfill \qed


The following result is used in the discussion of a Neumann boundary condition in Section \ref{ssecDNBC}.
\begin{lemma} \label{sumcvg1}
For $f \in L^{2}_{\rho^{(\beta \, , \, \alpha - \beta)}}(I)$ and $\mu_{i}$ and $c_{i}$ given by \eqref{defln} and \eqref{defuNs2}, 
respectively,
\be
\sum_{i = 0}^{\infty} \mu_{i} \, c_{i} \, G_{i+1}^{(\beta - 1 \, , \, \alpha - \beta - 1)}(0) \ < \ \infty \, .
\label{sume1}
\ee
\end{lemma}
\textbf{Proof}: Note that for $f \in L^{2}_{\rho^{(\beta \, , \, \alpha - \beta)}}(I)$,
\be
\infty \ > \ \| f \|_{L^{2}_{\rho^{(\beta \, , \, \alpha - \beta)}}}^{2} 
\ = \ \int_{0}^{1} \rho^{(\beta \, , \, \alpha - \beta)}(x) \, f(x)^{2} \, dx
\ = \ \sum_{i = 0}^{\infty} \frac{f_{i}^{2}}{ |\| G_{i}^{(\beta \, , \, \alpha - \beta)} |\|^{2} } \, .
\label{sume2}
\ee

From \eqref{defln} and \eqref{defuNs2} we have
\be
\mu_{i} \, c_{i}  \ = \ - \frac{\Gamma(i + \alpha)}{\Gamma(i + 1 + \alpha)} \, \frac{f_{i}}{ |\| G_{i}^{(\beta \, , \, \alpha - \beta)} |\|^{2} } \, .
\label{sume3}
\ee

Combining \eqref{eqG0}, \eqref{sume2} and \eqref{sume3},
\begin{align}
& \left| \sum_{i = 0}^{\infty} \mu_{i} \, c_{i} \, G_{i+1}^{(\beta - 1 \, , \, \alpha - \beta - 1)}(0) \right|
\ \le \ 
 \sum_{i = 0}^{\infty} \left| \, \mu_{i} \, c_{i} \, G_{i+1}^{(\beta - 1 \, , \, \alpha - \beta - 1)}(0) \right|  \nonumber \\
& \quad = \ \sum_{i = 0}^{\infty} \left|  \frac{ - \, \Gamma(i + \alpha)}{\Gamma(i + 1 + \alpha)} \,
 \frac{f_{i}}{ |\| G_{i}^{(\beta \, , \, \alpha - \beta)} |\|^{2} }  \, (-1)^{i+1} 
 \frac{ \Gamma(i + 1 + \alpha - \beta)}{ \Gamma(i + 2) \, \Gamma(\alpha - \beta)} \right|  \nonumber \\
& \quad \le \ 
\left( \sum_{i = 0}^{\infty}  \frac{f_{i}^{2}}{ |\| G_{i}^{(\beta \, , \, \alpha - \beta)} |\|^{2} } \right)^{1/2} \,
\left( \sum_{i = 0}^{\infty}  
\left( \frac{ \Gamma(i + \alpha)}{\Gamma(i + 1 + \alpha)}  \, 
\frac{ 1 }{ |\| G_{i}^{(\beta \, , \, \alpha - \beta)} |\|^{2} } \, 
 \frac{ \Gamma(i + 1 + \alpha - \beta)}{ \Gamma(i + 2) \, \Gamma(\alpha - \beta)}
  \right)^{2}  \, \right)^{1/2} \, .   
  \label{sume4}
\end{align}  

From \eqref{sume2} we have that the first term on the right hand side of \eqref{sume4} is bounded. Denote by
$S$ the second term on the right hand side of \eqref{sume4}. Using \eqref{spm22g},
\[
S^{2} \ = \ \frac{1}{\left(\Gamma(\alpha - \beta) \right)^{2}} \, 
\sum_{i = 0}^{\infty} \frac{1}{(i + \alpha)^{2}} \, 
\frac{(2 i \, + \, \alpha \, + \, 1)}{\Gamma(i + \beta + 1)} \, 
\frac{\Gamma(i + 1) \, \Gamma(i + \alpha + 1)}{\Gamma(i + \alpha - \beta + 1)} 
\left( \frac{\Gamma(i + 1 + \alpha - \beta)}{\Gamma(i + 2)} \right)^{2}  \, .
\]
Using Stirling's formula \eqref{eqStrf},
\[
\frac{\Gamma(i + 1)}{\Gamma(i + \beta + 1)} \ \sim \ (i + 1)^{-\beta} \, ,  \ \ 
\frac{\Gamma(i + \alpha + 1)}{\Gamma(i + \alpha - \beta + 1)} \ \sim \ (i + \alpha - \beta + 1)^{\beta} \, , \ \ 
\frac{\Gamma(i + 1 + \alpha - \beta)}{\Gamma(i + 2)}  \ \sim \ (i + 2)^{\alpha - \beta - 1} \, .
\]
Thus,
\[
    S^{2} \ \sim \ \frac{1}{\left(\Gamma(\alpha - \beta) \right)^{2}} \, \sum_{i = 0}^{\infty} i^{-1 \, + \, 2 ( \alpha - \beta - 1)} \ 
    < \ \infty \, , \ \  \mbox{ (as $\alpha - \beta - 1 \, < \, 0$) }
\]
from which the stated result then follows. \\
\mbox{ } \hfill \qed


\end{document}